\documentclass[lettersize,journal]{IEEEtran}

\usepackage{amsmath,amsfonts,amssymb}
\usepackage{mathtools}
\usepackage{bm}
\usepackage{breqn}
\usepackage{xcolor}
\usepackage{hyperref}
\usepackage{siunitx}

\usepackage{mathrsfs}

\usepackage{amsthm}
\theoremstyle{plain}%
\newtheorem{theorem}{Theorem}

\DeclareMathOperator\supp{supp}

\usepackage{rotating}
\usepackage{longtable}
\usepackage{booktabs}

\usepackage{algpseudocode}
\usepackage{algorithm}
\algrenewcommand{\algorithmicindent}{0.7em}  

\DeclareMathOperator{\RU}{RU}
\DeclareMathOperator{\RD}{RD}
\DeclareMathOperator{\UT}{UT}
\DeclareMathOperator{\DT}{DT}
\DeclareMathOperator{\RAE}{RAE}

\usepackage{pgfplots}
\pgfplotsset{compat=newest}
\usepgfplotslibrary{statistics}
\usepackage{xcolor}
\usepackage[export]{adjustbox} 
\usepackage{pgfplotstable} 
\usetikzlibrary{calc}

\makeatletter
\let\@msm@th@eqref\eqref
\renewcommand{\eqref}[1]{%
  \begingroup
  \leavevmode
  \color{blue}%
  \hypersetup{linkbordercolor=[named]{blue}}%
  \@msm@th@eqref{#1}%
  \endgroup
}
\makeatother

\usepackage{setspace}

\begin{document}

\title{Scenario Reduction for the Two-Stage Stochastic Unit Commitment Problem}

\author{Y. Werner, J. M. Morales, S. Pineda, L. Roald, and S. Wogrin \vspace{-0.4cm} %
\thanks{
Y. Werner and S. Wogrin are with the Institute of Electricity Economics and Energy Innovation at Graz University of Technology, 8020 Graz, Austria.\\
E-mails: $\{$yannick.werner, wogrin$\}$@tugraz.at.\\
J. M. Morales and S. Pineda are with the research group OASYS, University of Málaga, Málaga 29071, Spain.\\
E-mails: $\{$spineda,juan.morales$\}$@uma.es.\\
L. Roald is with the Department of Electrical and Computer Engineering, University of Wisconsin-Madison, Madison, WI 53706, USA.\\
E-mail: roald@wisc.edu.
}
}



\vspace{-1cm}
\IEEEaftertitletext{\vspace{-0.6\baselineskip}}
\maketitle

\begin{abstract}
The two-stage stochastic unit commitment problem has become an important tool to support decision-making under uncertainty in power systems. Representing the uncertainty by a large number of scenarios guarantees accurate results but challenges the solution process. One way to overcome this is by using scenario reduction methods, which aim at finding a distribution supported on fewer scenarios, but leading to similar optimal first-stage decisions. In this paper, we recap the classical scenario reduction theory based on the distance of probability distributions and the optimal mass transportation problem. We then review and compare various formulations of the underlying cost function of the latter used in the literature. Using the Forward Selection Algorithm, we show that a specific formulation of the cost function can be proven to select the best possible scenario from a given sample on the first draw with respect to the Relative Approximation Error. We demonstrate this result and compare the quality of the approximation as well as the computational performance of the different cost functions using a modified version of the IEEE RTS 24-Bus System. In many cases, we find that the optimal solution of the two-stage stochastic unit commitment problem with 200 scenarios can be approximated with around \qty{2}{\percent} error with only around 5 scenarios when using this cost function.
\end{abstract}

\begin{IEEEkeywords}
Sampling methods, Integer linear programming, Mathematical programming, Optimal scheduling, Forecast uncertainty
\end{IEEEkeywords}

\section{Introduction}
The increasing share of intermittent renewable power production challenges the economically optimal commitment of thermal generators \cite{Zheng2015}. The stochastic unit commitment problem (SUC) can support this decision-making process \cite{Takriti1996}. Frequently, the unknown distribution of uncertain renewable power production is approximated by a large enough set of scenarios. This has been recently motivated again by theoretical advancements on the sample average approximation \cite{Kleywegt2002} that provide mathematical guarantees for the asymptotic behavior of scenario-based approximations of continuous distributions (see also \cite{Shapiro2021}). At the same time, using a large set of scenarios increases the size of the optimization problem, potentially making it intractable or very difficult to solve, necessitating decomposition or other solution methods \cite{Zheng2015,Haaberg2019}.
 
Another way to overcome this problem is by applying scenario reduction techniques, which aim at reducing the number of scenarios needed to accurately represent the uncertainty in the stochastic problem without compromising its solution accuracy. Pioneering work by J. Dupačová et al.~\cite{Dupacova2003} related the efficiency of scenario reduction to the distance of probability distributions of the underlying random variables. They derive mathematical stability results and provide algorithms for scenario selection, most notably, the Forward Selection Algorithm, which we recap later on. This work was extended in~\cite{Heitsch2003} by improving the computational efficiency of the proposed algorithms and in~\cite{Heitsch2007} by considering different probability metrics. While earlier works largely addressed underlying stability and convergence properties (see also \cite{Pflug2001}, \cite{Roemisch2003} and \cite{Roemisch2007}), the applied distance measures only assess the similarity of input data and do not take the structure of the optimization problem into account, leading to slow convergence rates \cite{Henrion2018}. We refer to those cost functions as \textit{input-data-driven}. Several of these cost functions have also been successfully applied to the SUC problem. For example, a scenario reduction method inspired by importance sampling is developed in \cite{Papavasiliou2013}, and authors in~\cite{Dvorkin2014} present a comparison of different cost functions and selection techniques, such as k-means \cite{arthur2007k} and the backward scenario selection \cite{Heitsch2003}.

More recent literature has focused on incorporating information from the stochastic optimization problem into the scenario reduction process. We refer to those as \textit{optimization-problem-driven}. In~\cite{Morales2009}, the distance of scenarios is measured by comparing optimal objective function values of second-stage problems when fixing the first-stage decision variables to the optimal ones of the expected value problem. A related approach with fixed first-stage decisions is applied in~\cite{Feng2014} to the SUC problem with uncertain wind power production. Authors in \cite{Bruninx2016} show that for the SUC, this leads to risk-averse decision-making and propose using the difference in the objective function value of a single scenario deterministic two-stage problems instead. In~\cite{Bertsimas2023}, a symmetric cost function is proposed that takes the loss of decision quality into account when making first-stage decisions based on a specific scenario, while another one is realized. A similar approach is used in~\cite{Zhuang2025}, comparing different scenario reduction methods for the economic dispatch of active distribution networks and the SUC problem. 

Based on those works, we identify several gaps in the literature that challenge a rigorous comparison of different cost functions for the mass transportation problem. First, some papers proposing optimization-problem-driven cost functions or scenario reduction methods, such as \cite{Bertsimas2023} and \cite{Hewitt2021}, respectively, only compare their approach to input-data-driven methods. Second, papers that provide more in-depth comparisons, such as  \cite{Dvorkin2014} and \cite{Zhuang2025}, often vary distance measures and selection procedures at the same time. Finally, some papers proposing optimization-problem-driven distances that require calculating many small deterministic problems only consider convex \cite{Bertsimas2023} or small mixed-integer linear \cite{Zhuang2025} stochastic optimization problems.

Connecting to the classical literature on scenario reduction based on the distance of probability distributions and the optimal mass transportation problem, we try to bridge the aforementioned gaps by the following contributions:
\begin{itemize}
    \item We compare several optimization-problem-driven cost functions for the optimal mass transportation problem and evaluate them against a selected input-data-driven one.
    \item While varying the cost functions, we always use the Forward Selection Algorithm~\cite{Dupacova2003} to determine the reduced scenario set, allowing for a ceteris paribus comparison.
    \item For one of the cost functions, we develop and prove a mathematical theorem stating that the first scenario drawn with it is the best possible with respect to the approximation quality of the original distribution.
    \item Using a computationally challenging case study, we compare all cost functions on the SUC problem, evaluating the quality of the reduced scenario sets as well as the required computational expenses.
\end{itemize}

The remainder of this document is organized as follows. Section~\ref{sec:preliminaries} introduces a generic form of the two-stage SUC problem with uncertain wind power production. Afterwards, Section~\ref{sec:ScenarioReduction} recaps the classical theory on scenario reduction and the Forward Selection Algorithm~\cite{Dupacova2003}. Section~\ref{sec:CostFunctionMTP} introduces a selection of different cost functions found in the literature and states the optimality theorem for one of them. The performance of the cost functions applied to the two-stage SUC is demonstrated in Section~\ref{sec:results} and Section~\ref{sec:conclusions}.
\section{Preliminaries}
\label{sec:preliminaries}

In this paper, we apply scenario reduction methods to the two-stage SUC problem. For simplicity, we assume that the only uncertainty stems from intermittent renewable power production, which we represent by the random variable $\xi$. We further assume that $\xi$ can be adequately represented by a finite set of scenarios following the discrete distribution $\mathbb{P} = \sum_{i \in \mathcal{I}} p_i \delta_{\xi_i}$, where $p_i$ and $\xi_i \in \mathbb{R}^d$ are the probability and location of scenario $i$ for $i \in \mathcal{I} = \{1,\dots,n\}$ under distribution $\mathbb{P}$ and $\delta_{\xi_i}$ is the Dirac distribution allocating unit mass at $\xi_i$. Under these assumptions, we can represent the two-stage SUC problem with uncertain renewable power production, which we fully provide in Appendix~\ref{sec:app_SUC}, by the following generic two-stage stochastic program:
\begin{equation}\label{eq:GeneralTwoStageProblem}
    \min_{x \in \mathcal{X}} f(x) + \sum_{i \in \mathcal{I}} p_i G(x, \xi_i), 
\end{equation}
where $G(x, \xi_i)$ is given by:
\begin{subequations}\label{eq:GeneralSecondStage}
\begin{align}
    \min \quad & a^\top v \\
    \mathrm{s.t.} \quad & W v \geq h^{\xi_i} - T x, \label{eq:GeneralSecondStageConstraints}\\
    & v \geq 0. 
\end{align}
\end{subequations}
The goal of the two-stage stochastic program~\eqref {eq:GeneralTwoStageProblem} is choosing decision variables $x \in \mathcal{X}$, where $\mathcal{X}$ is defined by a set of mixed-integer linear constraints, to minimize first-stage cost $f(x)$ plus the expected second-stage cost taken with respect to the distribution of the random variable $\xi$, i.e., $\mathbb{E}_{\xi \sim \mathbb{P}}[G(x, \xi)] =  \sum_{i \in \mathcal{I}} p_i G(x, \xi_i)$, where the second stage problem~\eqref{eq:GeneralSecondStage} chooses continuous actions $v$ to minimize the linear operational cost $a^\top v$ subject to a set of linear constraints~\eqref{eq:GeneralSecondStageConstraints} that depend on the first-stage decisions $x$ and the realization $\xi_i$ of the random variable $\xi$. Notably, we take the following two common assumptions: (1) The cost function coefficients $a$, the recourse matrix $W$, and the technology matrix $T$ are independent of the uncertainty (note that the latter is generally not needed), and (2) the second-stage problem~\eqref{eq:GeneralSecondStage} is feasible for any realization $x \in \mathcal{X}$ (relatively complete recourse). For a discussion of those assumptions and the two-stage stochastic program in general, we refer the reader to~\cite{Birge2011} and \cite{Shapiro2021}.

In the context of the two-stage SUC problem, first-stage variables $x$ are the binary commitment decisions of dispatchable thermal units, which are subject to a set of technical constraints $\mathcal{X}$, e.g, minimum up- and downtimes, chosen to minimize a piece-wise linear cost function that takes, e.g., start-up and fixed costs into account. The second-stage problem minimizes the cost of operating the power system, ensuring technical feasibility, e.g., by restricting power production from generators of line power flows, after observing the realization of uncertain renewable power production $h^{\xi_i}$. 

We use $x^*(\mathbb{P})$ and $z^*(\mathbb{P})$ to denote the optimal first-stage decisions and objective function value, respectively, of Problem~\eqref{eq:GeneralTwoStageProblem} with respect to the distribution $\mathbb{P}$ as:
\begin{subequations}\label{eq:TwoStageOpt}
\begin{align}
    x^*(\mathbb{P}) = \arg&\min_{x \in \mathcal{X}} f(x) + \sum_{i \in \mathcal{I}} p_i G(x, \xi_i) \label{eq:TwoStageOptSol}\\
    z^*(\mathbb{P}) = z^*(x,\mathbb{P}) = &\min_{x \in \mathcal{X}} f(x) + \sum_{i \in \mathcal{I}} p_i G(x, \xi_i). \label{eq:TwoStageOptObj}
\end{align}
\end{subequations}

Additionally, we define for any given $\hat{x} \in \mathcal{X}$, $z^*(\hat{x},\mathbb{P})$, i.e., the optimal solution to Problem~\eqref{eq:GeneralTwoStageProblem} with  first-stage decisions fixed at $\hat{x}$, that is:
\begin{equation}
    z^*(\hat{x},\mathbb{P}) = f(\hat{x}) + \sum_{i \in \mathcal{I}} p_i G(\hat{x}, \xi_i).
\end{equation}

For the scenario reduction process later on, we further define a deterministic version of Problem~\eqref{eq:GeneralTwoStageProblem} considering a distribution supported on a single scenario $\xi_i$ only:
\begin{equation}\label{eq:GeneralSingleScenarioDP}
    \min_{x \in \mathcal{X}} z(x,\xi_i) = f(x) + G(x, \xi_i),
\end{equation}
where we denote the optimal solution to the first-stage decision considering scenario $\xi_i$ by:
\begin{equation}\label{eq:SingleScenarioOptSol}
    x^*(\xi_i) = \arg\min_{x \in \mathcal{X}} z(x,\xi_i).
\end{equation}

\section{A recap on scenario reduction}
\label{sec:ScenarioReduction}
The first thing we want to point out is that the real difficulty in stochastic programming in general is finding the optimal first-stage decisions $x^*(\mathbb{P})$. Once those are found, evaluating the individual, independent second-stage problems~\eqref{eq:GeneralSecondStage} to determine $z^*(\mathbb{P})$ is easy, as, for any given $x^*$ (in fact for any feasible $\hat{x} \in \mathcal{X}$), Equation~\eqref{eq:TwoStageOptObj} in combination with \eqref{eq:GeneralSingleScenarioDP} simplifies to:
\begin{equation}\label{eq:SecondStageEvaluationDistribution}
    z^*(x^*,\mathbb{P}) = \sum_{i \in \mathcal{I}} p_i (f(x^*) + G(x^*, \xi_i)) = \sum_{i \in \mathcal{I}} p_i z(x^*, \xi_i).
\end{equation}
Now, suppose there is another distribution $\mathbb{Q} = \sum_{j \in \mathcal{J}} q_j \delta_{\zeta_j}$, where $q_j$ and $\zeta_j \in \mathbb{R}^d$ are the probability and location of scenario $j$ for $j \in \mathcal{J} = \{1,\dots,m\}$ under distribution $\mathbb{Q}$ with $m < n$, and ideally $m \ll n$, as well as $x^*(\mathbb{Q}) \approx x^*(\mathbb{P})$.
Then we would not need to solve the optimization problem~\eqref{eq:GeneralTwoStageProblem} under distribution $\mathbb{P}$, which can be intractable or computationally challenging to solve, but under $\mathbb{Q}$, which is much easier. The goal of scenario reduction, which we discuss next, is to find such a distribution $\mathbb{Q}$.

\subsection{Distance of probability distributions}
\label{sec:ProbabilityDistributionDistance}
One way of assessing the similarity of probability distributions $\mathbb{P}$ and $\mathbb{Q}$ is by measuring their distance. In the mathematical programming literature, the distance $\mathfrak{D}$ is frequently quantified using the Monge-Kantorovich mass transportation problem (see e.g. \cite{Heitsch2003}):
\begin{equation}\label{eq:mass_transport}
\begin{aligned}
\mathfrak{D}(\mathbb{P}, \mathbb{Q}) = \min_{\pi \in \mathbb{R}_{+}^{n \times m}} \quad & \sum_{i=1}^{n} \sum_{j=1}^{m} \pi_{ij} c(\xi_i, \zeta_j) \\
\text{s.\,t.} \quad & \sum_{j=1}^{m} \pi_{ij} = p_i, \quad \forall i \in \mathcal{I}, \\
& \sum_{i=1}^{n} \pi_{ij} = q_j, \quad \forall j \in \mathcal{J},
\end{aligned}
\end{equation}
where $c(\xi_i, \zeta_j)$ is some function that measures the cost of moving the probability mass from scenario $\xi_i$ under distribution~$\mathbb{P}$ to $\zeta_j$ under distribution~$\mathbb{Q}$.
Let the support of distribution~$\mathbb{P}$ be defined as $\supp(\mathbb{P}) = \{\xi_i : i \in \mathcal{I} \}$ and by $\mathcal{P}(\Xi,m)$, for any set $\Xi \in \mathbb{R}^d$, the family of discrete distributions on $\Xi$ supported on $m$ points (scenarios). Then the goal of the \textit{discrete scenario reduction problem} is to find a distribution $\mathbb{Q} \in \mathcal{P}(\supp(\mathbb{P}),m)$, with $m < n$, and ideally $m \ll n$, such that both distributions $\mathbb{P}$ and $\mathbb{Q}$ are close with respect to some distance, e.g., as defined here in Equation~\eqref{eq:mass_transport} (cf. \cite{Rujeerapaiboon2018}).

\subsection{Forward Selection Algorithm}
\label{sec:ForwardSelectionAlgorithm}
In the seminal work on scenario reduction for stochastic programming using probability metrics, Dupačová et al. \cite{Dupacova2003} prove that in the case of the discrete scenario reduction problem, for any given set $\mathcal{J} \subset \mathcal{I}$, there is an analytical solution to Problem~\eqref{eq:mass_transport}, which only depends on $p_i$ and $c(\xi_i,\zeta_j)$ [\cite{Dupacova2003}, Theorem 3.1]:
\begin{equation}\label{eq:ProbDist_set}
    \mathfrak{D}(\mathbb{P}, \mathbb{Q}) = \mathfrak{D}(\mathbb{P}, \mathcal{J}) \coloneq \sum_{i \in \mathcal{I} \setminus \mathcal{J}} p_i \min_{j \in \mathcal{J}} c(\xi_i, \zeta_j),
\end{equation}
with
\begin{subequations}
\begin{align}
    q_j &= p_j + \sum_{i \in \mathcal{I}_j} p_i, \mathrm{where} \label{eq:ProbabiltyRedistribution}\\
    \mathcal{I}_j &\coloneq \{ i \in \mathcal{I}\,\setminus\,\mathcal{J} : j = \arg\min_{j' \in \mathcal{J}} ~ c(\xi_i, \zeta_{j'})\}.
\end{align}
\end{subequations}
Expression~\eqref{eq:ProbabiltyRedistribution} is referred to as the optimal probability redistribution rule. Intuitively, the redistribution of probabilities is such that the probability mass $p_i$ of all scenarios in $\mathcal{I}$ under distribution $\mathbb{P}$, which are closest to scenario $j$ with respect to $c(\xi_i, \zeta_j)$, are aggregated. The theorem implies that the goal of scenario reduction is actually finding the index set $J$, as the corresponding optimal probability distribution~$\mathbb{Q}$ follows straight away.

Based on that theorem, they propose the Forward Selection Algorithm \cite{Dupacova2003}, for which a pseudo-code is provided in Algorithm~\ref{alg:ForwardSelection}. We briefly recap the algorithm here for completeness, as it becomes relevant for the derivation of our theoretical result later on. To ease notation, we introduce an auxiliary set~$\mathcal{R}$ that captures all unselected scenarios from the original sample. Intuitively, the algorithm iteratively populates the set of scenarios $\mathcal{J}$ by a single scenario $j$ based on how the distance $\mathfrak{D}(\mathbb{P}, \mathcal{J})$, which is dependent on the choice of cost function $c$, would change if scenario $j$ were added to the reduced set of scenarios, until $m$ scenarios are selected in total. For clarity, we will use $\mathbb{P}_n$ and $\mathbb{Q}_m$ from here on to highlight the number of atoms the distributions are supported on.

While the Forward Selection Algorithm constitutes a greedy strategy, it is proven in~\cite{Dupacova2003} that it provides a globally optimal solution (selection) for $m=1$. We acknowledge that computational improvements to the Foward Selection Algorithm have been suggested in the literature, e.g., in~\cite{Heitsch2003}, but resort to the simpler version here as we find that for two-stage stochastic problems the computational time is generally negligible compared to solving the optimization problems.
\begin{algorithm}
\setstretch{1.15}
\caption{Forward Selection \cite{Dupacova2003}}
\begin{algorithmic}[1]
\Statex \textbf{inputs} $\mathbb{P}_n$, $m$, $c(\xi_i,\zeta_j)$
\State set $\mathcal{J} \gets \{\}$
\While{$|\mathcal{J}| < m$}
\State set $\mathcal{R} \gets \mathcal{I} \setminus \mathcal{J}$
\State $j = \arg\min_{j' \in \mathcal{R}} \sum_{i \in \mathcal{R} \setminus \{ j' \}} p_i \min_{j'' \in \mathcal{J} \cup \{j'\} } c(\xi_i,\zeta_{j''}) $
\State $\mathcal{J} \gets \mathcal{J} \cup \{j\}$
\EndWhile
\For{$j \in \mathcal{J}$}
    \State define $\mathcal{I}_j \coloneq \{ i \in \mathcal{I}\,\setminus\,\mathcal{J} : j = \arg\min_{j' \in \mathcal{J}} ~ c(\xi_i, \zeta_{j'}) \}$ 
    \State update $q_j = p_j + \sum_{i \in \mathcal{I}_j} p_i$
\EndFor
\Statex \textbf{return} $\mathbb{Q}_m$
\end{algorithmic}
\label{alg:ForwardSelection}
\end{algorithm}

\subsection{Evaluation of scenario reduction methods}
While the Forward Selection Algorithm~\eqref{alg:ForwardSelection} provides an efficient greedy heuristic to find a distribution~$\mathbb{Q}_m$, supported on $m$ scenarios, with a small distance~\eqref{eq:mass_transport}, its ability to accurately approximate the optimal decisions $x^*(\mathbb{P}_n)$ of Problem~\eqref{eq:GeneralSecondStage} hinges on the choice of the cost function $c(\xi_i,\zeta_j)$. Before discussing different formulations in Section~\ref{sec:CostFunctionMTP}, we briefly discuss how to evaluate the approximation quality of $\mathbb{Q}_m$ related to Problem~\eqref{eq:GeneralTwoStageProblem}.

For clarity, we define the optimal solution of the two-stage stochastic problem~\eqref{eq:GeneralTwoStageProblem} with respect to distribution $\mathbb{Q}_m$ in line with Equations~\eqref{eq:TwoStageOpt} as:
\begin{subequations}
\begin{align}
    x^*(\mathbb{Q}_m) &= \arg\min_{x \in \mathcal{X}} \sum_{j \in \mathcal{J}} q_j \cdot z(x,\zeta_j), \label{eq:TwoStageOptSolQ}  \\
    z^*(\mathbb{Q}_m), &= \min_{x \in \mathcal{X}} \sum_{j \in \mathcal{J}} q_j \cdot z(x,\zeta_j). \label{eq:TwoStageOptObjQ}
\end{align}
\end{subequations}

Building up on the definition in~\cite{Pflug2001} and \cite{kaut2003evaluation}, we then define the \textit{relative approximation error} $\RAE$ to quantify how well $\mathbb{Q}_m$ approximates $\mathbb{P}_n$ with respect to the two-stage stochastic Problem~\eqref{eq:GeneralTwoStageProblem}:
\begin{align}\label{eq:RAE}
    \RAE(\mathbb{P}_n,\mathbb{Q}_m) &= \frac{z^*(x^*(\mathbb{Q}_m),\mathbb{P}_n) - z^*(\mathbb{P}_n)}{z^*(\mathbb{P}_n)}, \\
    \mathrm{where}~z^*(x^*(\mathbb{Q}_m),\mathbb{P}_n) &= \sum_{i \in \mathcal{I}} p_i \cdot z(x^*(\mathbb{Q}_m),\xi_i). \nonumber
\end{align}
The $\RAE$ can be interpreted as the \textit{regret} of making decisions based on distribution $\mathbb{Q}_m$ while the true distribution is $\mathbb{P}_n$. 
We point out that this metric is only well defined when $z^*(\mathbb{P}_n) \neq 0$ and that it is always nonnegative for any $x \in \mathcal{X}$ by definition:
\begin{equation}
    z^*(\mathbb{P}_n) =  z^*(x^*(\mathbb{P}_n)), \mathbb{P}_n) \leq z^*(x, \mathbb{P}_n).
\end{equation}
This implies that there will always be a loss in decision quality when approximating distribution $\mathbb{P}_n$ with any distribution $\mathbb{Q}_m$ with respect to Problem~\eqref{eq:GeneralTwoStageProblem}, unless $x^*(\mathbb{Q}_m) = x^*(\mathbb{P}_n)$. Trivially, the $\RAE$ is zero when distributions $\mathbb{P}_n$ and $\mathbb{Q}_m$ are equal.

This metric has also recently been applied in~\cite{Zhuang2025} for the stochastic economic dispatch in active distribution networks and the SUC problem. We acknowledge that evaluating metric~\eqref{eq:RAE} is generally impractical, as the underlying assumption in scenario reduction is that it is not possible to solve the two-stage stochastic problem~\eqref{eq:GeneralTwoStageProblem} with distribution $\mathbb{P}_n$ in the first place. This fact is giving rise to the research on stability theorems \cite{Pflug2001, Dupacova2003, Roemisch2003}, focusing on establishing upper bounds to this error based on the distance of the probability distributions $\mathbb{P}_n$ and $\mathbb{Q}_n$. As our main focus here is to evaluate the impact of the choice of cost function $c$, however, we assume that we can actually solve Problem~\eqref{eq:GeneralTwoStageProblem} for comparison purposes.

\section{Cost functions for the transport problem}
\label{sec:CostFunctionMTP}
After recapping the scenario reduction process, we now want to take a look at various formulations for the cost function $c(\xi_i, \zeta_j)$ of the mass transportation problem~\eqref{eq:mass_transport} suggested in the literature. Note that in the context of stochastic programming, ideally we would like to define the cost function $c$ and distance $\mathfrak{D}$ directly with respect to the optimal first-stage decisions $x^*(\mathbb{P})$ and $x^*(\mathbb{Q})$ of Problem~\eqref{eq:GeneralTwoStageProblem}. However, for real-world applications, this is impractical, as we cannot solve those problems and therefore do not know the decisions $x^*(\mathbb{P})$ (or even $x^*(\mathbb{Q})$ for any distribution $\mathbb{Q}$). Therefore, we need to find a measure $c$ that provides a good way to approximate the loss in decision quality when moving from $x^*(\mathbb{P})$ to $x^*(\mathbb{Q})$.

We can generally separate the formulations proposed in the literature into two categories. First, those measures that are solely based on the realizations (scenarios) of the random variables $\xi_i$, i.e., the values of $h^{\xi_i}$ in Problem~\eqref{eq:GeneralSecondStage} discussed in this paper, which we refer to as \textit{input-data-driven}. And, secondly, those functions that take some information from the stochastic optimization problem~\eqref{eq:GeneralTwoStageProblem} or some approximation of it into account. We refer to those as \textit{optimization-problem-driven}. In the following, we introduce a selection of the most commonly discussed measures $c(\xi_i,\zeta_j)$.

\subsection{Input-data-driven}
In the context of input-data-driven measures, it has, for example, been proposed in~\cite{Rujeerapaiboon2018} to calculate the squared Euclidean norm of the scenario realizations:
\begin{equation}\label{eq:c_ID}
    c^{\mathrm{ID}}(\xi_i, \zeta_j) = \lVert \xi_i - \zeta_j \rVert_2^2.
\end{equation}
As e.g. shown in~\cite{Pflug2011}, the resulting distance $\mathfrak{D}$ is equal to the squared type-2 Wasserstein distance. We refer the interested reader to~\cite{Rujeerapaiboon2018}, where an in-depth description and analysis of the metric and its application to scenario reduction is provided. Here, we restrict ourselves to this distance in the class of input-data-driven measures, as it has also been used for comparison in other studies investigating optimization-problem-distances, such as \cite{Bertsimas2023}, and has been successfully applied to the two-stage SUC problem, e.g., in \cite{Dvorkin2014, Zhuang2025}.


\subsection{Optimization-problem-driven}
In the following, we introduce several optimization-problem-driven measures that, in contrast to the input-data-driven one in Equation~\eqref{eq:c_ID}, all rely on solving various versions and combinations of the single-scenario deterministic problem~\eqref{eq:GeneralSingleScenarioDP} and single-scenario second-stage problems~\eqref{eq:GeneralSecondStage} and therefore depend on the specific problem~\eqref{eq:GeneralTwoStageProblem} desired to solve.

\subsubsection{Morales et al. (2009)}
Authors in~\cite{Morales2009} suggest to measure similarity of scenarios based on the difference in second-stage costs when fixing the first-stage variables to those of the optimal solution of the expected value problem (EVP) \cite{Birge2011}:
\begin{equation}\label{eq:c_Mo}
    c^{\mathrm{Mo}}(\xi_i, \zeta_j) = | z(x^*(\bar{\xi}), \xi_i) - z(x^*(\bar{\xi}), \zeta_j) |,
\end{equation}
where $\bar{\xi} = \mathbb{E}_{\mathbb{P}_n}[\xi]$ is the expected realization of the random variable $\xi$ under distribution $\mathbb{P}_n$ and $x^*(\bar{\xi})$ is the optimal first-stage decision of the EVP.
When applied to the SUC, evaluating this measure requires solving a single MILP~\eqref{eq:GeneralSingleScenarioDP}, i.e., the EVP, and $n = |\mathcal{I}|$ second-stage LPs~\eqref{eq:GeneralSecondStage}.


\subsubsection{Bruninx et al. (2016)}
In~\cite{Bruninx2016}, it is pointed out that the measure $c^{\mathrm{Mo}}$ leads to a risk-averse selection of scenarios and, in the context of the SUC, tends to select scenarios with higher load shedding. To overcome this problem, they suggest using $c^{\mathrm{Br}}(\xi_i, \zeta_j)$, which considers the difference in optimal objective value of the single-scenario problems~\eqref{eq:GeneralSingleScenarioDP}: 
\begin{equation}\label{eq:c_Br}
    c^{\mathrm{Br}}(\xi_i, \zeta_j) = | z(x^*(\xi_i), \xi_i) - z(x^*(\zeta_j), \zeta_j) |.
\end{equation}
Evaluating this measure requires solving $n$ MILP problems.

\subsubsection{Bertsimas et al. (2023)}
Aligned with the stability theory presented in the classical literature on scenario reduction using input-data-driven distances, e.g., \cite{Dupacova2003,Heitsch2003}, it is suggested in~\cite{Bertsimas2023} to use distance:
\begin{align}\label{eq:c_Be}
    2 \cdot c^{\mathrm{Be}}(\xi_i, \zeta_j) = z&(x^*(\zeta_j), \xi_i) - z(x^*(\xi_i),\xi_i) + \nonumber\\
    z&(x^*(\xi_i), \zeta_j) - z(x^*(\zeta_j),\zeta_j),
\end{align}
which considers the symmetric loss in decision quality when making first-stage decisions based on scenario $\xi_i$ but then observing $\zeta_j$ and vice versa. In order to evaluate measure $c^{\mathrm{Be}}$, one needs to solve $n$ MILPs~\eqref{eq:GeneralSingleScenarioDP} and $n^2 - n$ second-stage LPs~\eqref{eq:GeneralSecondStage}.


\subsubsection{Proposed}
In this paper, we propose using measure $c^{\mathrm{Pr}}$ to calculate~\eqref{eq:mass_transport}:
\begin{equation}\label{eq:c_Pr}
    c^{\mathrm{Pr}}(\xi_i, \zeta_j) = z(x^*(\zeta_j), \xi_i) - z(x^*(\xi_i),\xi_i),
\end{equation}
which measures the nonsymmetric regret when determining the optimal first-stage decisions based on scenario $\zeta_j$, while the actual realization is $\xi_i$. Like distance $c^{\mathrm{Be}}$, it requires solving $n$ MILPs~\eqref{eq:GeneralSingleScenarioDP} and $n^2 - n$ second-stage LPs~\eqref{eq:GeneralSecondStage}.

Before further discussing the properties of this measure, we acknowledge that the same cost function was first proposed in~\cite{Hewitt2021}, referred to as opportunity cost, and in \cite{Keutchayan2023}, referred to as prediction error, in the context of nonhierarchical scenario clustering. While \cite{Hewitt2021} focuses on deriving a bounding strategy for the clustering algorithm, the main goal in \cite{Keutchayan2023} is to find good cluster representatives. In stark contrast to their works, we apply the measure $c^{\mathrm{Pr}}$ to the classical scenario reduction framework based on the distance of probability distributions in~\eqref{eq:mass_transport} and the greedy Forward Selection Algorithm~\eqref{alg:ForwardSelection}. This allows us to establish the following theorem, which is our main theoretical contribution:
\begin{theorem}\label{th:Theorem1}
    Let $\mathcal{Q}_1 = \{ \mathbb{Q}_1^j : j \in \mathcal{I} \}$ denote the family of probability distributions obtained by drawing a single scenario from distribution $\mathbb{P}_n$ and $j^*$ the first scenario drawn using the Forward Selection Algorithm~\ref{alg:ForwardSelection} and cost function $c^{\mathrm{Pr}}$. Then the distribution $\mathbb{Q}_1^{j^*}$ minimizes the $\RAE$ for $m = 1$, such that $\mathbb{Q}_1^{j^*} = \arg\min_{\mathbb{Q}_1 \in \mathcal{Q}_1} \mathrm{RAE}(\mathbb{P}_n, \mathbb{Q}_1)$.
\end{theorem}
\textit{Proof:} See Appendix~\ref{sec:app_proof}. 
\vspace{1mm}

\noindent The intuitive interpretation of Theorem~\eqref{th:Theorem1} is that the Forward Selection Algorithm~\eqref{alg:ForwardSelection}, when using cost function~\eqref{eq:c_Pr}, selects the best possible scenario from distribution $\mathbb{P}_n$ with respect to the $\RAE$ on the very first draw. An intuitive explanation for this is that when using the cost function $c^{\mathrm{Pr}}$, the Forward Selection Algorithm~\ref{alg:ForwardSelection} actually compares the $\RAE$ associated with drawing each individual scenario first, and selects exactly the one scenario that minimizes it. Even though this result can only be established for the first scenario drawn, we show in Section~\ref{sec:results} that for the SUC and the case studies used here, this first scenario constitutes a good starting point for finding a suitable distribution $\mathbb{Q}_m$.

\subsection{Summary}
We summarize the required optimization problems that must be solved to evaluate the cost functions $\mathcal{C} = \{ \mathrm{ID},\mathrm{Mo},\mathrm{Br}, \mathrm{Be}, \mathrm{Pr} \}$ in Table~\ref{tab:OptProbRequirments}.
\begin{table}[ht]
\centering
\caption{Comparison of optimization problems needed to solve to evaluate cost function $c$.}
\label{tab:OptProbRequirments}
\begin{tabular}{lccccc}
\toprule
 & ID & Mo & Br & Be & Pr \\
\midrule
MILPs & $0$ & $1$ & $n$ & $n$ & $n$ \\
LPs & $0$ & $n$ & $0$ & $n^2 - n$ & $n^2 - n$ \\
\bottomrule
\end{tabular}
\end{table}
In Section~\ref{sec:CaseStudy_24Bus}, we provide a detailed analysis of the computational expense of calculating those problems for the two-stage SUC problem and the case studies used here.

Algorithm~\ref{alg:ScenarioReduction} provides a pseudo-code for the whole scenario reduction process, including solving the reduced two-stage stochastic problems and evaluating the solution. Again, the goal of scenario reduction is to find a distribution $\mathbb{Q}_m$ based on an original distribution $\mathbb{P}_n$, with $m < n$ and ideally $m \ll n$, that provides a good approximation $x^*(\mathbb{Q}_m) \approx x^*(\mathbb{P}_n)$ with respect to the stochastic optimization problem~\eqref{eq:GeneralTwoStageProblem}, evaluated based on the $\RAE$ in Equation~\eqref{eq:RAE}.

\begin{algorithm}[]
\setstretch{1.15}
\caption{Scenario Reduction}
\begin{algorithmic}[1]
\State \textbf{inputs} $\mathbb{P}_n$, $m$
\Statex \textit{Step 0: Initialization}
\State choose $c \in \{ c^{\mathrm{ID}}, c^{\mathrm{Mo}}, c^{\mathrm{Br}}, c^{\mathrm{Be}},   c^{\mathrm{Pr}}\}$
\Statex \textit{Step 1:} \textit{Calculate cost matrix}
\For{$i \in \mathcal{I}$, $j \in \mathcal{I}$}
    \State calculate $c(\xi_i, \zeta_j)$ after solving problems in Table~\ref{tab:OptProbRequirments}
\EndFor
\Statex \textit{Step 2:} \textit{Determine reduced distribution $\mathbb{Q}_m$}
\State run Forward Selection Algorithm~\ref{alg:ForwardSelection} with ($\mathbb{P}_n$, $m$, $c(\xi_i, \zeta_j)$), receive $\mathbb{Q}_m$
\Statex \textit{Step 3:} \textit{Solve two-stage stochastic problem}
\State solve $x^*(\mathbb{Q}_m) \in \arg\min_{x \in \mathcal{X}} ~ z(x,\mathbb{Q}_m)$
\Statex \textit{Step 4:} \textit{Evaluate first-stage decisions}
\State calculate $\RAE(\mathbb{P}_n,\mathbb{Q}_m)$
\State \textbf{return} $\mathbb{Q}_m$, $\RAE(\mathbb{P}_n,\mathbb{Q}_m)$
\end{algorithmic}
\label{alg:ScenarioReduction}
\end{algorithm}

\section{Numerical results}
\label{sec:results}
In the following, we analyze the numerical results for applying the Scenario Reduction Algorithm~\ref{alg:ScenarioReduction} to the SUC problem provided in Appendix~\ref{sec:app_SUC}, while varying the cost functions $c$ introduced in Section~\ref{sec:CostFunctionMTP}. Everything is implemented in the Julia Programming language v1.11.5 \cite{bezanson2017julia} using JuMP v1.23.3 \cite{Lubin2023} and all optimization models are solved with Gurobi v.11.0.2 \cite{gurobi}. The code and the case study data can be found on GitHub \cite{github_SRSUC}. All simulations are executed on an Intel Core i9-13900 Processor using 24 cores clocking at a base rate of \qty{2.00}{\giga\hertz} and up to a maximum of \qty{5.60}{\giga\hertz} in turbo mode.


%
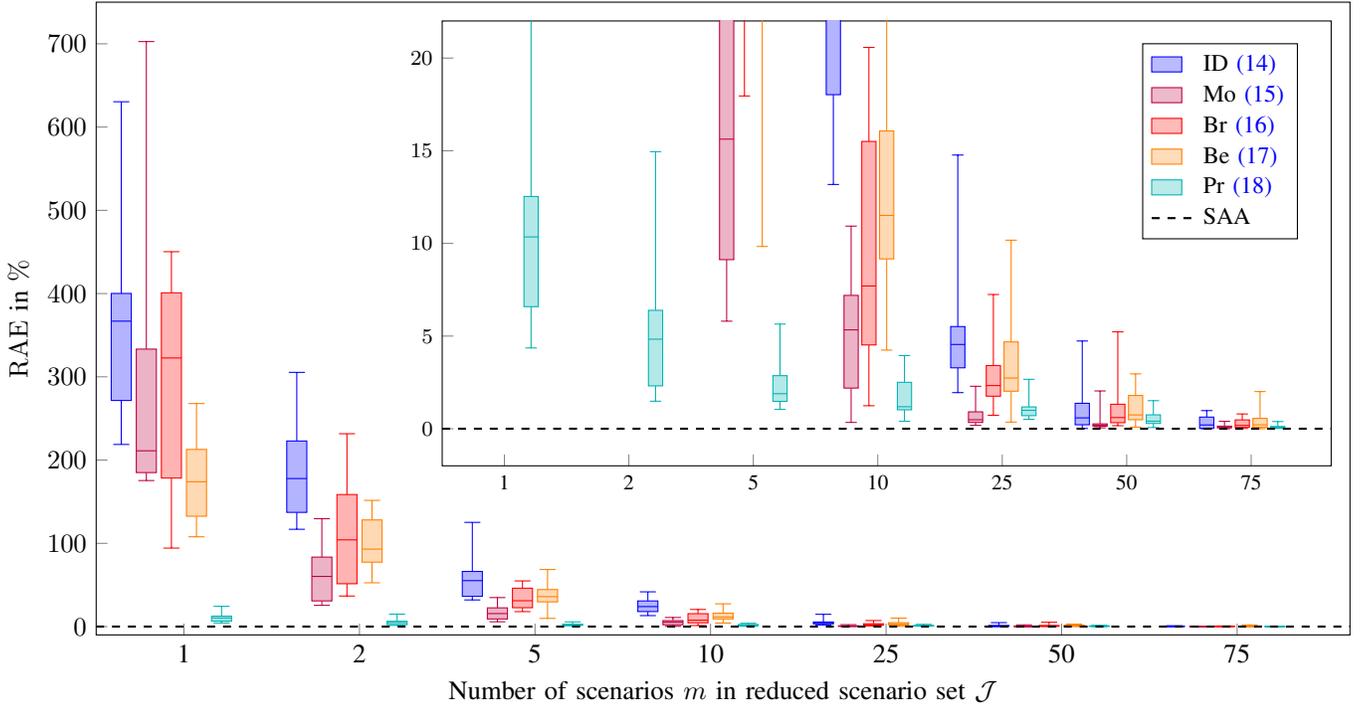
\begin{figure*}

\pgfplotstableread[col sep=comma]{plots/data/rae_boxplot_data.csv}\datatable

\pgfplotstablegetrowsof{\datatable}
\pgfmathsetmacro{\numrows}{\pgfplotsretval-1}

\colorlet{colIDDSR}{blue}
\colorlet{colMorales}{purple}
\colorlet{colBruninx}{red}
\colorlet{colBertsimas}{orange}
\colorlet{colBertsimasUnstable}{teal!60!cyan}

\newcommand{\methodlist}{IDDSR,Morales,Bruninx,Bertsimas,BertsimasUnstable}
\newcommand{\methodlabels}{IDDSR,Morales,Bruninx,Bertsimas,Proposed}
\newcommand{\clusterlist}{1,2,5,10,25,50,75}
\newcommand{\nummethods}{5}
\newcommand{\groupwidth}{7}

\begin{tikzpicture}
\begin{axis}[
    name=mainplot,
    boxplot/draw direction=y,
    xmin=0, xmax=50,
    ymin=-10, ymax=750,
    xtick={3.5, 10.5, 17.5, 24.5, 31.5, 38.5, 45.5},
    xticklabels={1, 2, 5, 10, 25, 50, 75},
    xlabel={Number of scenarios $m$ in reduced scenario set $\mathcal{J}$},
    ylabel={$\mathrm{RAE}$ in \unit{\percent}},
    width=18.25cm,
    height=10cm,
    ymajorgrids=false,
    grid style={gray!30},
    xtick pos=bottom,
    ytick pos=left,
]

\pgfplotsinvokeforeach{0,...,\numrows}{
    \pgfplotstablegetelem{#1}{method}\of\datatable
    \let\currentmethod\pgfplotsretval
    \pgfplotstablegetelem{#1}{cluster}\of\datatable
    \let\currentcluster\pgfplotsretval
    \pgfplotstablegetelem{#1}{min}\of\datatable
    \let\currentmin\pgfplotsretval
    \pgfplotstablegetelem{#1}{q1}\of\datatable
    \let\currentqone\pgfplotsretval
    \pgfplotstablegetelem{#1}{median}\of\datatable
    \let\currentmedian\pgfplotsretval
    \pgfplotstablegetelem{#1}{q3}\of\datatable
    \let\currentqthree\pgfplotsretval
    \pgfplotstablegetelem{#1}{max}\of\datatable
    \let\currentmax\pgfplotsretval
    
    \ifnum\pdfstrcmp{\currentmethod}{IDDSR}=0
        \def\methodidx{0}
        \def\methodcolor{colIDDSR}
    \fi
    \ifnum\pdfstrcmp{\currentmethod}{Morales}=0
        \def\methodidx{1}
        \def\methodcolor{colMorales}
    \fi
    \ifnum\pdfstrcmp{\currentmethod}{Bruninx}=0
        \def\methodidx{2}
        \def\methodcolor{colBruninx}
    \fi
    \ifnum\pdfstrcmp{\currentmethod}{Bertsimas}=0
        \def\methodidx{3}
        \def\methodcolor{colBertsimas}
    \fi
    \ifnum\pdfstrcmp{\currentmethod}{BertsimasUnstable}=0
        \def\methodidx{4}
        \def\methodcolor{colBertsimasUnstable}
    \fi
    
    \ifnum\currentcluster=1
        \def\clusteridx{0}
    \fi
    \ifnum\currentcluster=2
        \def\clusteridx{1}
    \fi
    \ifnum\currentcluster=5
        \def\clusteridx{2}
    \fi
    \ifnum\currentcluster=10
        \def\clusteridx{3}
    \fi
    \ifnum\currentcluster=25
        \def\clusteridx{4}
    \fi
    \ifnum\currentcluster=50
        \def\clusteridx{5}
    \fi
    \ifnum\currentcluster=75
        \def\clusteridx{6}
    \fi
    
    \pgfmathsetmacro{\xpos}{1 + \clusteridx*7 + \methodidx}
    
    \edef\temp{\noexpand\addplot+[
        boxplot prepared={
            lower whisker=\currentmin,
            lower quartile=\currentqone,
            median=\currentmedian,
            upper quartile=\currentqthree,
            upper whisker=\currentmax,
            draw position=\xpos,
        },
        fill=\methodcolor!30,
        draw=\methodcolor,
        solid,
        forget plot,
    ] coordinates {};}
    \temp
}

\addplot[black, thick, dashed, domain=0:50, samples=2, mark=none, forget plot] {0};

\end{axis}

\begin{axis}[
    name=insetplot,
    at={(mainplot.north east)},
    anchor=north east,
    xshift=-0.25cm,
    yshift=-0.25cm,
    boxplot/draw direction=y,
    xmin=0, xmax=50,
    ymin=-2, ymax=22,
    xtick={3.5, 10.5, 17.5, 24.5, 31.5, 38.5, 45.5},
    xticklabels={1, 2, 5, 10, 25, 50, 75},
    xlabel={},
    ylabel={},
    width=13.4cm,
    height=7.5cm,
    ymajorgrids=false,
    grid style={gray!30},
    axis background/.style={fill=white},
    tick label style={font=\footnotesize},
    xtick pos=bottom,
    ytick pos=left,
    legend style={
        at={(0.875,0.95)},
        anchor=north,
        legend columns=1,
        cells={anchor=west},
        font=\small,
        fill=white,
        fill opacity=0.9,
        draw opacity=1,
        text opacity=1,
    },
    legend image code/.code={
        \draw[#1, fill] (0cm,-0.1cm) rectangle (0.4cm,0.1cm);
    },
]

\addlegendimage{fill=colIDDSR!30, draw=colIDDSR}
\addlegendentry{ID~\eqref{eq:c_ID}}
\addlegendimage{fill=colMorales!30, draw=colMorales}
\addlegendentry{Mo~\eqref{eq:c_Mo}}
\addlegendimage{fill=colBruninx!30, draw=colBruninx}
\addlegendentry{Br~\eqref{eq:c_Br}}
\addlegendimage{fill=colBertsimas!30, draw=colBertsimas}
\addlegendentry{Be~\eqref{eq:c_Be}}
\addlegendimage{fill=colBertsimasUnstable!30, draw=colBertsimasUnstable}
\addlegendentry{Pr~\eqref{eq:c_Pr}}
\addlegendimage{line legend, black, thick, dashed}
\addlegendentry{SAA}

\pgfplotsinvokeforeach{0,...,\numrows}{
    \pgfplotstablegetelem{#1}{method}\of\datatable
    \let\currentmethod\pgfplotsretval
    \pgfplotstablegetelem{#1}{cluster}\of\datatable
    \let\currentcluster\pgfplotsretval
    \pgfplotstablegetelem{#1}{min}\of\datatable
    \let\currentmin\pgfplotsretval
    \pgfplotstablegetelem{#1}{q1}\of\datatable
    \let\currentqone\pgfplotsretval
    \pgfplotstablegetelem{#1}{median}\of\datatable
    \let\currentmedian\pgfplotsretval
    \pgfplotstablegetelem{#1}{q3}\of\datatable
    \let\currentqthree\pgfplotsretval
    \pgfplotstablegetelem{#1}{max}\of\datatable
    \let\currentmax\pgfplotsretval
    
    \ifnum\pdfstrcmp{\currentmethod}{IDDSR}=0
        \def\methodidx{0}
        \def\methodcolor{colIDDSR}
    \fi
    \ifnum\pdfstrcmp{\currentmethod}{Morales}=0
        \def\methodidx{1}
        \def\methodcolor{colMorales}
    \fi
    \ifnum\pdfstrcmp{\currentmethod}{Bruninx}=0
        \def\methodidx{2}
        \def\methodcolor{colBruninx}
    \fi
    \ifnum\pdfstrcmp{\currentmethod}{Bertsimas}=0
        \def\methodidx{3}
        \def\methodcolor{colBertsimas}
    \fi
    \ifnum\pdfstrcmp{\currentmethod}{BertsimasUnstable}=0
        \def\methodidx{4}
        \def\methodcolor{colBertsimasUnstable}
    \fi
    
    \ifnum\currentcluster=1
        \def\clusteridx{0}
    \fi
    \ifnum\currentcluster=2
        \def\clusteridx{1}
    \fi
    \ifnum\currentcluster=5
        \def\clusteridx{2}
    \fi
    \ifnum\currentcluster=10
        \def\clusteridx{3}
    \fi
    \ifnum\currentcluster=25
        \def\clusteridx{4}
    \fi
    \ifnum\currentcluster=50
        \def\clusteridx{5}
    \fi
    \ifnum\currentcluster=75
        \def\clusteridx{6}
    \fi
    
    \pgfmathsetmacro{\xpos}{1 + \clusteridx*7 + \methodidx}
    
    \pgfmathsetmacro{\clippedmin}{max(\currentmin, -2)}
    \pgfmathsetmacro{\clippedmax}{min(\currentmax, 100)}
    \pgfmathsetmacro{\clippedqone}{max(min(\currentqone, 100), -2)}
    \pgfmathsetmacro{\clippedmedian}{max(min(\currentmedian, 100), -2)}
    \pgfmathsetmacro{\clippedqthree}{max(min(\currentqthree, 100), -2)}
    
    \edef\temp{\noexpand\addplot+[
        boxplot prepared={
            lower whisker=\clippedmin,
            lower quartile=\clippedqone,
            median=\clippedmedian,
            upper quartile=\clippedqthree,
            upper whisker=\clippedmax,
            draw position=\xpos,
        },
        fill=\methodcolor!30,
        draw=\methodcolor,
        solid,
        forget plot,
    ] coordinates {};}
    \temp
}

\addplot[black, thick, dashed, domain=0:50, samples=2, mark=none, forget plot] {0};

\end{axis}

\end{tikzpicture}

    \caption{$\RAE$ for the cost functions $c$ presented in Section~\ref{sec:CostFunctionMTP} depending on the size of the reduced scenario set $\mathcal{J}$, using the Scenario Reduction Algorithm~\ref{alg:ScenarioReduction}. Box plots represent variations over different initial sample distributions $\mathbb{P}_n$, each of size $n = 200$, drawn from the same probabilistic forecast~\cite{pinson2013} using the modified IEEE RTS 24-Bus System described in Section~\ref{sec:CaseStudy_24Bus}.}
    \label{fig:IEEE24_quality_draws}
\end{figure*}

\subsection{Case study: IEEE RTS 24-Bus System}
\label{sec:CaseStudy_24Bus}

The first case study is a modified version of the IEEE RTS 24-bus System presented in \cite{Ordoudis2016}, where we use the original 32 dispatchable generators of the original data set \cite{Subcommittee1979}, which we assume must be all committed. We also add a generator-specific fixed cost term $C^{\mathrm{fix}}$ and modify the linear production costs $C^{\mathrm{L}}$ based on the version presented in the PGLib~\cite{Babaeinejadsarookolaee2019}. Finally, we slightly perturb all cost coefficients to ensure that each generator has a unique cost function and reduce the potential for multiple optimal solutions.
The power system hosts six wind power generators at nodes \{3,5,7,16,21,23\} with an installed capacity of \SI{350}{\mega\watt} each. We assume that the uncertain wind power production can be approximated by a probabilistic forecast with up to 10000 equiprobable scenarios based on the data in \cite{pinson2013}, which is based on real measurements of wind farms in Denmark.
We apply uniform load shedding cost $C^{\mathrm{Sh}}$ of \qty{1500}{\$\per\mega\watt\hour} and assume that wind power can be curtailed at no cost.
All values are converted to per unit such that the optimal objective function value of Problem~\eqref{eq:GeneralTwoStageProblem} is in the second order of magnitude. We use Gurobi's standard solver settings if not indicated otherwise.

\subsection{Quality of approximation}

To assess the impact of using the different cost functions $c$ discussed in Section~\ref{sec:CostFunctionMTP}, we take 11 different random draws of size $|\mathcal{I}| = 200$ from the available probabilistic forecast and run the Scenario Reduction Algorithm~\ref{alg:ScenarioReduction}. Figure~\ref{fig:IEEE24_quality_draws} shows the distribution of the $\RAE$ over the draws for different sizes $m$ of the reduced set of scenarios $J$. The dashed line denotes the $\RAE$ when solving the two-stage SUC for distribution $\mathbb{P}_n$, which is equal to the sample average approximation (SAA).

The first thing that we would like to point out is that for the case where only a single scenario is drawn from the original sample, such that $m=1$, all cost functions discussed in Section~\ref{sec:CostFunctionMTP}, except $c^{\mathrm{Pr}}$, find quite inaccurate reduced sets of scenarios with low approximation quality and large $\RAE$. In particular, depending on the sample drawn, they largely vary in approximation quality. For the cost function $c^{\mathrm{Pr}}$, on the other hand, we can observe two things: (1) The numerical simulations confirm the result in Theorem~\ref{th:Theorem1}, i.e., that always the best possible scenario in terms of $\RAE$ is selected, leading to a comparably accurate approximation of the first stage decisions for $m=1$, and (2) that because it is guaranteed to select this scenario independent of the original distribution $\mathbb{P}_n$ drawn, it delivers much more robust results when varying the sample. In fact, it is often possible to achieve a $\RAE$ of around \qty{10}{\percent} with only a single scenario. While Theorem~\ref{th:Theorem1} only considers the case of $m=1$, we can see that the distributions for a higher number of scenarios delivered by the cost function $c^{\mathrm{Pr}}$ show good quality and convergence. As shown in the inset plot of Figure~\ref{fig:IEEE24_quality_draws}, all cost functions start delivering distributions with good approximation quality from 25 scenarios onward for most of the draws. Depending on the desired level of accuracy, however, we find that when using $c^{\mathrm{Pr}}$, the optimal first-stage decisions of the two-stage SUC problem related to the original distribution $\mathbb{P}_n$ can be approximated with an $\RAE$ of around \qty{2}{\percent} with only 5 scenarios. Notably, we also find that $c^{\mathrm{Mo}}$ delivers more accurate results than cost functions $c^{\mathrm{Br}}$ and $c^{\mathrm{Be}}$, while being much easier to evaluate, which is what we analyze next.

\subsection{Computational expenses}

After analyzing the quality of the distributions produced by the different cost functions $c$ in Figure~\eqref{fig:IEEE24_quality_draws}, we now take a look at the computational expenses required to evaluate the various cost functions. Table~\ref{tab:computation_times} shows the total Gurobi work units to solve the necessary optimization problems shown in Table~\ref{tab:OptProbRequirments} for the 24-Bus System used here.
\begin{table}
\centering
\caption{Mean and standard deviation (in brackets) of Gurobi's work units required for solving the optimization problems in \textit{Step~1} of Algorithm~\ref{alg:ScenarioReduction}, taken over all draws.}
\label{tab:computation_times}
\begin{tabular}{lcc}
\toprule
Cost func. $c$ & MILPs & LPs \\
\midrule
ID~\eqref{eq:c_ID} & --- & --- \\
Mo~\eqref{eq:c_Mo} & 3.49 (1.75) & 3.57 (0.15) \\
Br~\eqref{eq:c_Br} & 4361.14 (615.24) & --- \\
Be~\eqref{eq:c_Be} & 4361.14 (615.24) & 720.46 (27.52) \\
Pr~\eqref{eq:c_Pr} & 4361.14 (615.24) & 720.46 (27.52) \\
\bottomrule
\end{tabular}
\end{table}
Since $c^{\mathrm{ID}}$ is evaluated on input data only, no optimization problems need to be solved. For the cost function $c^{\mathrm{Mo}}$, solving the EVP takes around 3.49 work units on average, while it takes around the same work to solve all the $|\mathcal{I}| = 200$ second-stage operational LPs. For cost functions $c^{\mathrm{Br}}$, $c^{\mathrm{Be}}$, and $c^{\mathrm{Pr}}$, it is necessary to solve 200 MILPs, which takes around 4360 work units in total. For the latter two, solving the additional around $200^2$ LPs takes another 720 work units in total. There are a couple of points we would like to highlight here. First, solving the individual MILP problems is computationally much more expensive than solving a comparably larger number of LP problems, and exhibits large variations depending on the initial sample. Therefore, secondly, while we find that the cost function $c^{\mathrm{Pr}}$ yields good reduced sets of scenarios, a tradeoff may arise between accuracy and computational time when considering two-stage stochastic MILP problems. Authors in~\cite{Keutchayan2023} argue that one may solve relaxed MILPs to evaluate $c$. We find this to be a very sensitive heuristic that may largely depend on the problem instance. Finally, taking the quality of the approximation into account, we would like to point out that the cost function $c^{\mathrm{Mo}}$ may be a good choice when the overall available time is very limited, and some error in the approximation can be tolerated, especially compared to the cost function $c^{\mathrm{ID}}$, which is solely input-data-driven.
\section{Summary and conclusions}
\label{sec:conclusions}
In this paper, we apply scenario reduction to a two-stage SUC problem with uncertain wind power production to find a distribution that is supported on fewer atoms than the originally given sample without changing the corresponding optimal unit commitment decisions too much. To find the reduced set of scenarios, we use the Forward Selection Algorithm~\cite{Dupacova2003} that iteratively evaluates the distance of probability distributions if a scenario $j$ would be selected based on the optimal mass transportation problem and a given cost function.

We provide an in-depth analysis of various cost functions suggested in the literature and prove a mathematical theorem for the one we propose to use, which states that it selects the best possible scenario on the first draw with respect to the $\RAE$. Applying the two-stage SUC to a modified version of IEEE RTS 24-Bus System, we compare the quality of approximation and computational expense for each of the cost functions. We find that the proposed cost function often achieves a high approximation quality of around \qty{2}{\percent} with 5 scenarios only. However, for two-stage stochastic MILPs, such as the SUC problem, it may take substantial time to evaluate the cost function, as it requires solving a series of smaller MILPs and LPs before. 

Future work should therefore focus on reducing the computational expense of evaluating those problems by developing efficient warm-starting strategies, such as using surrogate models. Furthermore, for cases where the second-stage problems are convex, one may investigate if and how duality information can be leveraged in the cost function. Finally, it should be analyzed how to best derive lower bounds to the two-stage stochastic problem for a given reduced set of scenarios, to enable the evaluation of the approximation quality without solving the original problem.

\appendices

\setcounter{equation}{0}
\renewcommand\theequation{A.\arabic{equation}}

\section{Stochastic unit commitment problem}
\label{sec:app_SUC}
This appendix presents the formulation of the two-stage SUC problem, following~\cite{Carrion2006} and \cite{Blanco2017}. We use lower-case letters for variables and upper-case letters for parameters, if not indicated otherwise.
Let $n, m \in \mathcal{N}$, $g \in \mathcal{G}$, $j \in \mathcal{J}$, and $t \in \mathcal{T} = \{1,2,...,T\}$ denote the set of nodes (buses), dispatchable generators, wind power generators, and time periods, respectively. Furthermore, let $\omega \in \Omega$ denote the set of wind power scenarios with associated wind power production $\widetilde{W}_{j,t,\omega}$.

We use $T^{U^0}_g$ to denote the number of time steps that generator $g$ has already been operating in its state $U^0_g$ at $t = 1$ and $\mathcal{T}^{0}_g$ for the number of time steps that a generator $g$ has to remain in its initial state to account for minimum up and down times $\UT$ and $\DT$, respectively, as:
\begin{equation}
    T^{0}_g =
    \begin{cases}
        \min (T, \UT_g - T^{U^0}_g) & \text{if}~U^0_g = 1, \\
        \min (T, \DT_g - T^{U^0}_g) & \text{if}~U^0_g = 0.
    \end{cases}
\end{equation}

For example, if a generator $g'$ has been offline for 8 time steps prior to $t = 1$, such that $U^0_g = 0$ and $T^{U^0}_{g'} = 8$, and its minimum down time is $\DT_{g'} = 12$ then $T^{0}_g = 4$ and $\mathcal{T}^{0}_{g'} = \{1,2,3,4\}$.

The goal of the unit commitment problem is to optimally choose commitment decisions $u_{g,t}$ (on/off), start-ups $y_{g,t}$. and shut-downs $z_{g,t}$ for each dispatchable generator $g$ and time step $t$. Here, the objective functionof the first-stage problem is given by:
\begin{equation}
    \sum_{t \in \mathcal{T}} \sum_{g \in \mathcal{G}} \left( C_g^{\mathrm{SU}} y_{g,t} + C_g^{\mathrm{fix}} u_{g,t} \right), \label{eq:FirstStage_Objective}
\end{equation}
where $C_g^{\mathrm{SU}}$ and $C_g^{\mathrm{fix}}$ refer to the start-up and fix costs, respectively.

The feasible region of the first-stage problem is given by:
\begin{subequations}\label{eq:FirstStage_feasible_region}
\begin{align}
& y_{g,t} - z_{g,t} = u_{g,t} - u_{g,t-1}, & & \forall g , \forall t \in \{2,...,T\}, \label{eq:FirstStage_state_change} \\
& y_{g,t} - z_{g,t} = u_{g,t} - U_g^0, & & \forall g, t = 1, \label{eq:FirstStage_state_change_init} \\
& y_{g,t} + z_{g,t} \leq 1, & & \forall g , \forall t \in \mathcal{T}, \label{eq:FirstStage_state_change_single} \\
& u_{g,t} = U_g^0 , & & \forall g , \forall t \in \mathcal{T}^0_g, \label{eq:FirstStage_state_initial}  \\
& \sum_{\tau = t - \mathrm{UT}_g + 1}^t y_{g,t} \leq u_{g,t}, & & \forall g, \forall t \in \mathcal{T} \setminus \mathcal{T}^0_g, \label{eq:FirstStage_min_up_time} \\
& \sum_{\tau = t - \mathrm{DT}_g + 1}^t z_{g,t} \leq 1 - u_{g,t}, & & \forall g, \forall t \in \mathcal{T} \setminus \mathcal{T}^0_g, \label{eq:FirstStage_min_down_time}
\end{align}
\end{subequations}
where $x = \{  (u_{g,t}, y_{g,t}, z_{g,t}), \forall g, \forall t\}$. Constraints~\eqref{eq:FirstStage_state_change}--\eqref{eq:FirstStage_state_change_init} relate start-ups and shut-downs to on/off states, Constraints~\eqref{eq:FirstStage_state_change_single} forbid simultaneous start-up and shut-down, Constraints~\eqref{eq:FirstStage_state_initial} enforce the initial state to an exogenously given one $U_g^0$, and Constraints~\eqref{eq:FirstStage_min_up_time}--\eqref{eq:FirstStage_min_down_time} define minimum up and down times, respectively.

The objective function for the second-stage problem of the SUC problem for each scenario $\omega \in \Omega$ is given by:
\begin{equation}\label{eq:SecondStage_objective}
   \sum_{t \in \mathcal{T}} \left( \sum_{g \in \mathcal{G}} C_g^{\mathrm{L}} p_{g,t,\omega} + \sum_{d \in \mathcal{D}} C^{\mathrm{Sh}} l^{\mathrm{Sh}}_{d,t,\omega} + \sum_{n \in \mathcal{N}} C_n^{\mathrm{Sl}} l^{\mathrm{Sl}}_{n,t,\omega} \right),
\end{equation}
where$p_{g,t,\omega}$ is the power production of dispatchable generator $g$ in time step $t$ and scenario $\omega$ at cost $C_g^{\mathrm{L}}$, $l^{\mathrm{Sh}}_{d,t,\omega}$ is the load shedding of demand $d$ at uniform cost $C^{\mathrm{Sh}}$, and $l^{\mathrm{Sl}}_{n,t,\omega}$ is a nodal excess slack demand that ensures feasibility for every first-stage unit commitment decision at cost $C_d^{\mathrm{Sl}}$, which we assume to be equal to the uniform load shedding cost $C^{\mathrm{Sh}}$.

The feasible region for every scenario $\omega \in \Omega$ is given by:
\begin{subequations}\label{eq:SecondStage_constraints}
\begin{align}
& P_g^{\min} u_{g,t} \leq p_{g,t,\omega} \leq P_g^{\max} u_{g,t}, & & \forall g, \forall t \in \mathcal{T}, \label{eq:SecondStage_prod_limit} \\
& -\RD_g \leq p_{g,t,\omega} - p_{g,t-1,\omega}, & & \forall g, \forall t \in \{2,...,T\}, \label{eq:SecondStage_ramp_limit_down} \\
& p_{g,t,\omega} - p_{g,t-1,\omega} \leq \RU_g, & & \forall g, \forall t \in \{2,...,T\}, \label{eq:SecondStage_ramp_limit_up} \\
& (P_g^{0} + \RD_g) u_{g,t} \leq p_{g,t,\omega}, & & \forall g, t = 1, \label{eq:SecondStage_ramp_limit_init_down} \\
& p_{g,t,\omega} \leq (P_g^{0} + \RU_g) u_{g,t}, & & \forall g, t = 1, \label{eq:SecondStage_ramp_limit_init_up} \\
& 0 \leq l^{\mathrm{Sh}}_{d,t,\omega} \leq L_{d,t}, & & \forall d, \forall t \in \mathcal{T}, \label{eq:SecondStage_shed_limit} \\
& 0 \leq w^{\mathrm{Sp}}_{j,t,\omega} \leq W_{j,t,\omega}, & & \forall j, \forall t \in \mathcal{T}, \label{eq:SecondStage_spill_limit} \\
& f_{n,m,t,\omega} = B_{n,m} (\theta_{n,t,\omega} - \theta_{m,t,\omega}), & & \forall n, \forall m \in \mathcal{N}_n, \forall t \in \mathcal{T}, \label{eq:SecondStage_flow_def} \\
& -\mathrm{F}_{n,m}^{\max} \leq f_{n,m,t,\omega} \leq \mathrm{F}_{n,m}^{\max}, & & \forall n, \forall m \in \mathcal{N}_n, \forall t \in \mathcal{T},\label{eq:SecondStage_line_limit} \\
& \sum_{g \in \mathcal{G}_n} p_{g,t,\omega}  - \sum_{m \in \mathcal{N}_n} f_{n,m,t,\omega} - & &  \nonumber \\
& l^{\mathrm{Sl}}_{n,t,\omega} + \sum_{j \in \mathcal{J}_n} (\widetilde{W}_{j,t,\omega} - w^{\mathrm{Sp}}_{j,t,\omega}) - & & \nonumber \\
& \sum_{d \in \mathcal{D}_n} (L_{d,t} - l^{\mathrm{Sh}}_{d,t,\omega})  = 0, & & \forall n, \forall t \in \mathcal{T}, \label{eq:SecondStage_nodal_bal}
\end{align}
\end{subequations}
where Constraints~\eqref{eq:SecondStage_prod_limit} limit power production of dispatchable generator $g$ in time step $t$ and scenario $\omega$ to be between the minimum and maximum technical limits $P_g^{\min}$ and $P_g^{\max}$, respectively. Constraints~\eqref{eq:SecondStage_ramp_limit_down}--\eqref{eq:SecondStage_ramp_limit_init_up} restrict their ramping capabilities to between a maximum ramp-down and ramp-up limit $\RD_g$ and $\RU_g$, respectively. Load shedding $l^{\mathrm{Sh}}_{d,t,\omega}$ of demand $d$ is restricted to be below the exogenously given inflexible load $L_{d,t}$ by Constraints~\eqref{eq:SecondStage_shed_limit}. Constraints~\eqref{eq:SecondStage_spill_limit} limit the amount of wind power spillage $w^{\mathrm{Sp}}_{j,t,\omega}$ to be below the respective forecast in scenario $\omega$, $\widetilde{W}_{j,t,\omega}$. We use a loss less DC power flow approximation to define line flows $f_{n,m,t,\omega}$ based on line susceptance $B_{n,m}$ and nodal voltage angles $\theta_{n,t,\omega}$ through Constraints~\eqref{eq:SecondStage_flow_def}. Symmetric line flow limits are enforced by Constraints~\eqref{eq:SecondStage_line_limit} and nodal power balance is ensured by Constraints~\eqref{eq:SecondStage_nodal_bal}.

The second-stage operational problem of the SUC for every scenario $\omega \in \Omega$ is then given by:
\begin{equation}\label{eq:SecondStage}
    \min_{} \eqref{eq:SecondStage_objective}, \mathrm{s.t.} \eqref{eq:SecondStage_constraints}.
\end{equation}
Note that depending on whether the second-stage problem is solved without the first-stage problem, then $u_{g,t}$ will be a parameter, not a variable.

The full two-stage SUC problem is then given by:
\begin{equation}\label{eq:TwoStageUnitCommitment}
\min_{x \in \eqref{eq:FirstStage_feasible_region}} \eqref{eq:FirstStage_Objective} + \mathbb{E}_{\omega}[\eqref{eq:SecondStage}] = \min_{x \in \eqref{eq:FirstStage_feasible_region}} \eqref{eq:FirstStage_Objective} + \sum_{\omega \in \Omega} \pi_\omega \cdot \eqref{eq:SecondStage}, 
\end{equation}
where the equality follows from the assumption that $\omega$ follows a finite discrete distribution where $\pi_\omega$ denotes the probability of scenario $\omega$.

\section{Proof of Theorem 1}
\label{sec:app_proof}

We first note that term $z^*(\mathbb{P}_n)$ in the definition of the $\RAE$ in~\eqref{eq:RAE} only depends on the distribution $\mathbb{P}_n$ and not on any set of scenarios $\mathcal{J}$ drawn from it. Consequently, the solution $x^*(\mathbb{Q}_m)$, defined by Equation~\eqref{eq:TwoStageOptSolQ}, that minimizes the first term of the nominator also minimizes the $\RAE$. Secondly, there exists only one unique probability distribution carried on a single atom: the Dirac measure that places unit probability mass on that same atom, i.e.,  $\mathbb{Q}_1^j = \delta_{\zeta_j}, \forall j \in \mathcal{I}$. Now observe that in this case whenever the distribution of $\xi$ is supported solely on a single scenario $\zeta_j$, such that $\xi \sim \mathbb{Q}_1^j = \delta_{\zeta_j}$ and $p_j = 1$, the general two-stage stochastic problem~\eqref{eq:GeneralTwoStageProblem} and the single-scenario problem~\eqref{eq:GeneralSingleScenarioDP} are identical. Therefore, their optimal solutions $x^*(\mathbb{Q}_1^j)$ and $x^*(\zeta_j)$,  defined by Equations~\eqref{eq:TwoStageOptSol} and \eqref{eq:SingleScenarioOptSol}, respectively, must also be identical, i.e., $x^*(\mathbb{Q}_1^j) = x^*(\zeta_j)$. Based on those points and assuming\footnote{While we find this assumption to be naturally fulfilled in most problems, it can always be ensured by adding a large enough positive constant to the objective function.} that $z^*(\mathbb{P}_n) > 0$ we derive:
\begin{align}
    \mathbb{Q}_1^{j^*} &= \arg\min_{\mathbb{Q}_1^j \in \hat{\mathcal{Q}}_1} \RAE(\mathbb{P}_n, \mathbb{Q}_1^j) \nonumber \\
    &= \arg\min_{\mathbb{Q}_1^j \in \hat{\mathcal{Q}}_1} \frac{z(x^*(\mathbb{Q}_1^j),\mathbb{P}_n) - z^*(\mathbb{P}_n)}{z^*(\mathbb{P}_n)} \nonumber \\
    &= \arg\min_{\mathbb{Q}_1^j \in \hat{\mathcal{Q}}_1} z(x^*(\mathbb{Q}_1^j),\mathbb{P}_n) \nonumber\\
    &= \arg\min_{j \in \mathcal{I}} z(x^*(\zeta_j),\mathbb{P}_n) \nonumber \\
    &= \arg\min_{j \in \mathcal{I}} \sum_{i \in \mathcal{I}} p_i z(x^*(\zeta_j),\xi_i), \label{eq:RAEReformulationZeta}
\end{align}
where the last equivalence follows from Equation~\eqref{eq:SecondStageEvaluationDistribution}.

Now taking a look at the first scenario drawn from $\mathcal{I}$ using the FS Algorithm~\ref{alg:ForwardSelection}, such that $\mathcal{J} = \{\}$ and $\mathcal{R} = \mathcal{I}$, Line~4 can be rewritten as:
\begin{align}
    j &= \arg\min_{j' \in \mathcal{R}} \sum_{i \in \mathcal{R} \setminus \{ j' \}} p_i \min_{j'' \in \mathcal{J} \cup \{j'\} } c(\xi_i,\zeta_{j''}) \nonumber\\
    &= \arg\min_{j' \in \mathcal{I}} \sum_{i \in \mathcal{I} \setminus \{ j' \}} p_i \min_{j'' \in \{\} \cup \{j'\} } c(\xi_i,\zeta_{j''}) \nonumber \\
    &= \arg\min_{j' \in \mathcal{I}} \sum_{i \in \mathcal{I} \setminus \{ j' \}} p_i c(\xi_i,\zeta_{j'}) \nonumber \\
    &= \arg\min_{j' \in \mathcal{I}} \sum_{i \in \mathcal{I}} p_i c(\xi_i,\zeta_{j'}), \label{eq:FSReformulation}
\end{align}
where the last equivalence follows from $c(\xi_{j'},\zeta_{j'}) = 0$.

Finally, we observe that for the measure $c^{\mathrm{Pr}}$, the second term, $z(x(\xi_i),\xi_i))$, is independent of scenario $j$, and, similar to the $\RAE$, is therefore identical for every scenario $j$. Now, replacing the alias $j'$ with $j$ and setting $c(\xi_i,\zeta_{j}) =  c^{\mathrm{Pr}}(\xi_i,\zeta_{j})$ in Equation~\eqref{eq:FSReformulation} we obtain:
\begin{align}
    j^* &= \arg\min_{j \in \mathcal{I}} \sum_{i \in \mathcal{I}} p_i c^{\mathrm{Pr}}(\xi_i,\zeta_{j}) \nonumber \\
    &= \arg\min_{j \in \mathcal{I}} \sum_{i \in \mathcal{I}} p_i (z(x^*(\zeta_{j}), \xi_i) - z(x^*(\xi_i),\xi_i)) \nonumber \\
    &= \arg\min_{j \in \mathcal{I}} \sum_{i \in \mathcal{I}} p_i z(x^*(\zeta_{j}), \xi_i), \label{eq:FSOptimalityFirstScenario}
\end{align}
which is equivalent to the minimization of the $\RAE$ in Equation~\eqref{eq:RAEReformulationZeta} and completes the proof. $\qedsymbol$

\bibliographystyle{IEEEtran}
\bibliography{paper/SR_for_UC}

@Article{Bertsimas2023,
  author    = {Bertsimas, Dimitris and Mundru, Nishanth},
  journal   = {Operations Research},
  title     = {Optimization-Based Scenario Reduction for Data-Driven Two-Stage Stochastic Optimization},
  year      = {2023},
  issn      = {1526-5463},
  month     = jul,
  number    = {4},
  pages     = {1343--1361},
  volume    = {71},
  doi       = {10.1287/opre.2022.2265},
  publisher = {Institute for Operations Research and the Management Sciences (INFORMS)},
}

@Book{Birge2011,
  author    = {Birge, John R. and Louveaux, François},
  publisher = {Springer New York},
  title     = {Introduction to Stochastic Programming},
  year      = {2011},
  isbn      = {9781461402374},
  doi       = {10.1007/978-1-4614-0237-4},
  issn      = {1431-8598},
  journal   = {Springer Series in Operations Research and Financial Engineering},
}

@InProceedings{Bruninx2016,
  author    = {Bruninx, K. and Delarue, E.},
  booktitle = {2016 IEEE International Energy Conference (ENERGYCON)},
  title     = {Scenario reduction techniques and solution stability for stochastic unit commitment problems},
  year      = {2016},
  month     = apr,
  pages     = {1--7},
  publisher = {IEEE},
  doi       = {10.1109/energycon.2016.7514074},
}

@Article{Dupacova2003,
  author    = {Jitka Dupačová and Nicole Gröwe-Kuska and Werner Römisch},
  journal   = {Mathematical Programming},
  title     = {Scenario reduction in stochastic programming},
  year      = {2003},
  issn      = {1436-4646},
  month     = mar,
  number    = {3},
  pages     = {493--511},
  volume    = {95},
  doi       = {10.1007/s10107-002-0331-0},
  publisher = {Springer Science and Business Media LLC},
}

@InProceedings{Feng2014,
  author    = {Feng, Yonghan and Ryan, Sarah M.},
  booktitle = {2014 IEEE PES General Meeting | Conference \& Exposition},
  title     = {Scenario reduction for stochastic unit commitment with wind penetration},
  year      = {2014},
  month     = jul,
  pages     = {1--5},
  publisher = {IEEE},
  doi       = {10.1109/pesgm.2014.6939138},
  ranking   = {rank1},
}

@Article{Haaberg2019,
  author    = {Håberg, Martin},
  journal   = {International Journal of Electrical Power \& Energy Systems},
  title     = {Fundamentals and recent developments in stochastic unit commitment},
  year      = {2019},
  issn      = {0142-0615},
  month     = jul,
  pages     = {38--48},
  volume    = {109},
  doi       = {10.1016/j.ijepes.2019.01.037},
  publisher = {Elsevier BV},
}

@Article{Heitsch2003,
  author    = {Heitsch, Holger and Römisch, Werner},
  journal   = {Computational Optimization and Applications},
  title     = {Scenario Reduction Algorithms in Stochastic Programming},
  year      = {2003},
  issn      = {0926-6003},
  number    = {2/3},
  pages     = {187--206},
  volume    = {24},
  doi       = {10.1023/a:1021805924152},
  publisher = {Springer Science and Business Media LLC},
}

@Article{Heitsch2007,
  author    = {Heitsch, Holger and Römisch, Werner},
  journal   = {Operations Research Letters},
  title     = {A note on scenario reduction for two-stage stochastic programs},
  year      = {2007},
  issn      = {0167-6377},
  month     = nov,
  number    = {6},
  pages     = {731--738},
  volume    = {35},
  doi       = {10.1016/j.orl.2006.12.008},
  publisher = {Elsevier BV},
}

@Article{Henrion2018,
  author    = {Henrion, R. and Römisch, W.},
  journal   = {Mathematical Programming},
  title     = {Problem-based optimal scenario generation and reduction in stochastic programming},
  year      = {2018},
  issn      = {1436-4646},
  month     = oct,
  number    = {1},
  pages     = {183--205},
  volume    = {191},
  doi       = {10.1007/s10107-018-1337-6},
  publisher = {Springer Science and Business Media LLC},
}

@Article{Hewitt2021,
  author    = {Hewitt, Mike and Ortmann, Janosch and Rei, Walter},
  journal   = {Annals of Operations Research},
  title     = {Decision-based scenario clustering for decision-making under uncertainty},
  year      = {2021},
  issn      = {1572-9338},
  month     = jan,
  number    = {2},
  pages     = {747--771},
  volume    = {315},
  doi       = {10.1007/s10479-020-03843-x},
  publisher = {Springer Science and Business Media LLC},
}

@Article{Keutchayan2023,
  author    = {Keutchayan, Julien and Ortmann, Janosch and Rei, Walter},
  journal   = {Computational Management Science},
  title     = {Problem-driven scenario clustering in stochastic optimization},
  year      = {2023},
  issn      = {1619-6988},
  month     = mar,
  number    = {1},
  volume    = {20},
  doi       = {10.1007/s10287-023-00446-2},
  publisher = {Springer Science and Business Media LLC},
}

@Article{Kleywegt2002,
  author    = {Kleywegt, Anton J. and Shapiro, Alexander and Homem-de-Mello, Tito},
  journal   = {SIAM Journal on Optimization},
  title     = {The Sample Average Approximation Method for Stochastic Discrete Optimization},
  year      = {2002},
  issn      = {1095-7189},
  month     = jan,
  number    = {2},
  pages     = {479--502},
  volume    = {12},
  doi       = {10.1137/s1052623499363220},
  publisher = {Society for Industrial & Applied Mathematics (SIAM)},
}

@Article{Morales2009,
  author    = {Morales, J.M. and Pineda, S. and Conejo, A.J. and Carrion, M.},
  journal   = {IEEE Transactions on Power Systems},
  title     = {Scenario Reduction for Futures Market Trading in Electricity Markets},
  year      = {2009},
  issn      = {1558-0679},
  month     = may,
  number    = {2},
  pages     = {878--888},
  volume    = {24},
  doi       = {10.1109/tpwrs.2009.2016072},
  publisher = {Institute of Electrical and Electronics Engineers (IEEE)},
}

@Book{Ordoudis2016,
  author    = {Christos Ordoudis and Pierre Pinson and {Morales Gonz{\'a}lez}, {Juan Miguel} and Marco Zugno},
  publisher = {Technical University of Denmark},
  title     = {An Updated Version of the IEEE RTS 24-Bus System for Electricity Market and Power System Operation Studies.},
  year      = {2016},
  language  = {English},
}

@InBook{Pflug2011,
  author    = {Pflug, Georg Ch. and Pichler, Alois},
  pages     = {343--387},
  publisher = {Springer New York},
  title     = {Approximations for Probability Distributions and Stochastic Optimization Problems},
  year      = {2011},
  isbn      = {9781441995865},
  booktitle = {Stochastic Optimization Methods in Finance and Energy},
  doi       = {10.1007/978-1-4419-9586-5_15},
  issn      = {0884-8289},
}

@Article{Pinson2013,
  author    = {Pinson, Pierre},
  journal   = {Statist. Sci.},
  title     = {Wind Energy: Forecasting Challenges for Its Operational Management},
  year      = {2013},
  month     = {11},
  number    = {4},
  pages     = {564--585},
  volume    = {28},
  doi       = {10.1214/13-STS445},
  fjournal  = {Statistical Science},
  publisher = {The Institute of Mathematical Statistics},
  url       = {https://doi.org/10.1214/13-STS445},
}

@InBook{Roemisch2003,
  author    = {Römisch, Werner},
  pages     = {483--554},
  publisher = {Elsevier},
  title     = {Stability of Stochastic Programming Problems},
  year      = {2003},
  isbn      = {9780444508546},
  booktitle = {Stochastic Programming},
  doi       = {10.1016/s0927-0507(03)10008-4},
  issn      = {0927-0507},
}

@Article{Roemisch2007,
  author    = {Römisch, W. and Wets, R. J.-B.},
  journal   = {SIAM Journal on Optimization},
  title     = {Stability of $\varepsilon$-approximate Solutions to Convex Stochastic Programs},
  year      = {2007},
  issn      = {1095-7189},
  month     = jan,
  number    = {3},
  pages     = {961--979},
  volume    = {18},
  doi       = {10.1137/060657716},
  publisher = {Society for Industrial & Applied Mathematics (SIAM)},
}

@Article{Rujeerapaiboon2018,
  author    = {Rujeerapaiboon, Napat and Schindler, Kilian and Kuhn, Daniel and Wiesemann, Wolfram},
  journal   = {Mathematical Programming},
  title     = {Scenario reduction revisited: fundamental limits and guarantees},
  year      = {2018},
  issn      = {1436-4646},
  month     = apr,
  number    = {1},
  pages     = {207--242},
  volume    = {191},
  doi       = {10.1007/s10107-018-1269-1},
  publisher = {Springer Science and Business Media LLC},
}

@Book{Shapiro2021,
  author    = {Shapiro, Alexander},
  editor    = {Darinka Dentcheva and Andrzej P. Ruszczyński},
  publisher = {Society for Industrial and Applied Mathematics},
  title     = {Lectures on stochastic programming},
  year      = {2021},
  address   = {Philadelphia},
  edition   = {Third edition},
  isbn      = {9781611976595},
  number    = {28},
  series    = {MOS-SIAM series on optimization},
  pagetotal = {1525},
  ppn_gvk   = {1774673193},
  subtitle  = {Modeling and theory},
}

@Article{Subcommittee1979,
  author    = {Subcommittee, Probability},
  journal   = {IEEE Transactions on Power Apparatus and Systems},
  title     = {{IEEE} Reliability Test System},
  year      = {1979},
  issn      = {0018-9510},
  month     = nov,
  number    = {6},
  pages     = {2047--2054},
  volume    = {PAS-98},
  doi       = {10.1109/tpas.1979.319398},
  publisher = {Institute of Electrical and Electronics Engineers (IEEE)},
}

@Article{Takriti1996,
  author    = {Takriti, S. and Birge, J.R. and Long, E.},
  journal   = {IEEE Transactions on Power Systems},
  title     = {A stochastic model for the unit commitment problem},
  year      = {1996},
  issn      = {0885-8950},
  number    = {3},
  pages     = {1497--1508},
  volume    = {11},
  doi       = {10.1109/59.535691},
  publisher = {Institute of Electrical and Electronics Engineers (IEEE)},
}

@Article{Zheng2015,
  author    = {Zheng, Qipeng P. and Wang, Jianhui and Liu, Andrew L.},
  journal   = {IEEE Transactions on Power Systems},
  title     = {Stochastic Optimization for Unit Commitment—A Review},
  year      = {2015},
  issn      = {1558-0679},
  month     = jul,
  number    = {4},
  pages     = {1913--1924},
  volume    = {30},
  doi       = {10.1109/tpwrs.2014.2355204},
  publisher = {Institute of Electrical and Electronics Engineers (IEEE)},
}

@Article{Zhuang2025,
  author    = {Zhuang, Yingrui and Cheng, Lin and Qi, Ning and Almassalkhi, Mads R. and Liu, Feng},
  journal   = {IEEE Transactions on Power Systems},
  title     = {Problem-Driven Scenario Reduction Framework for Power System Stochastic Operation},
  year      = {2025},
  issn      = {1558-0679},
  month     = jul,
  number    = {4},
  pages     = {3232--3246},
  volume    = {40},
  doi       = {10.1109/tpwrs.2024.3523220},
  publisher = {Institute of Electrical and Electronics Engineers (IEEE)},
}

@Book{kaut2003evaluation,
  author    = {Kaut, Michal and Wallace, Stein W},
  publisher = {Humboldt-Universit{\"a}t zu Berlin, Berlin},
  title     = {Evaluation of scenario-generation methods for stochastic programming},
  year      = {2003},
}

@InProceedings{Dvorkin2014,
  author    = {Dvorkin, Yury and Wang, Yishen and Pandzic, Hrvoje and Kirschen, Daniel},
  booktitle = {2014 IEEE PES General Meeting | Conference \& Exposition},
  title     = {Comparison of scenario reduction techniques for the stochastic unit commitment},
  year      = {2014},
  month     = jul,
  pages     = {1--5},
  publisher = {IEEE},
  doi       = {10.1109/pesgm.2014.6939042},
}

@Article{Pflug2001,
  author    = {Pflug, G.Ch.},
  journal   = {Mathematical Programming},
  title     = {Scenario tree generation for multiperiod financial optimization by optimal discretization},
  year      = {2001},
  issn      = {0025-5610},
  month     = jan,
  number    = {2},
  pages     = {251--271},
  volume    = {89},
  doi       = {10.1007/pl00011398},
  publisher = {Springer Science and Business Media LLC},
}

@Article{Carrion2006,
  author    = {Carrion, M. and Arroyo, J.M.},
  journal   = {IEEE Transactions on Power Systems},
  title     = {A Computationally Efficient Mixed-Integer Linear Formulation for the Thermal Unit Commitment Problem},
  year      = {2006},
  issn      = {0885-8950},
  month     = aug,
  number    = {3},
  pages     = {1371--1378},
  volume    = {21},
  doi       = {10.1109/tpwrs.2006.876672},
  publisher = {Institute of Electrical and Electronics Engineers (IEEE)},
}

@Article{Blanco2017,
  author    = {Blanco, Ignacio and Morales, Juan M.},
  journal   = {IEEE Transactions on Power Systems},
  title     = {An Efficient Robust Solution to the Two-Stage Stochastic Unit Commitment Problem},
  year      = {2017},
  issn      = {1558-0679},
  month     = nov,
  number    = {6},
  pages     = {4477--4488},
  volume    = {32},
  doi       = {10.1109/tpwrs.2017.2683263},
  publisher = {Institute of Electrical and Electronics Engineers (IEEE)},
}

@Article{Papavasiliou2013,
  author    = {Papavasiliou, Anthony and Oren, Shmuel S.},
  journal   = {Operations Research},
  title     = {Multiarea Stochastic Unit Commitment for High Wind Penetration in a Transmission Constrained Network},
  year      = {2013},
  issn      = {1526-5463},
  month     = jun,
  number    = {3},
  pages     = {578--592},
  volume    = {61},
  doi       = {10.1287/opre.2013.1174},
  publisher = {Institute for Operations Research and the Management Sciences (INFORMS)},
}

@InProceedings{arthur2007k,
  author    = {Arthur, David and Vassilvitskii, Sergei},
  booktitle = {Proceedings of the Eighteenth Annual ACM-SIAM Symposium on Discrete Algorithms (SODA)},
  title     = {k-means++: The Advantages of Careful Seeding},
  year      = {2007},
  pages     = {1027--1035},
}

@Misc{Babaeinejadsarookolaee2019,
  author    = {Babaeinejadsarookolaee, Sogol and Birchfield, Adam and Christie, Richard D. and Coffrin, Carleton and DeMarco, Christopher and Diao, Ruisheng and Ferris, Michael and Fliscounakis, Stephane and Greene, Scott and Huang, Renke and Josz, Cedric and Korab, Roman and Lesieutre, Bernard and Maeght, Jean and Mak, Terrence W. K. and Molzahn, Daniel K. and Overbye, Thomas J. and Panciatici, Patrick and Park, Byungkwon and Snodgrass, Jonathan and Tbaileh, Ahmad and Van Hentenryck, Pascal and Zimmerman, Ray},
  title     = {The Power Grid Library for Benchmarking AC Optimal Power Flow Algorithms},
  year      = {2019},
  copyright = {Creative Commons Attribution 4.0 International},
  doi       = {10.48550/ARXIV.1908.02788},
  publisher = {arXiv},
}

@Misc{github_SRSUC,
  author    = {Werner, Yannick},
  title     = {Scenario Reduction for the Two-Stage Stochastic Unit Commitment Problem},
  year      = {2025},
  publisher = {GitHub},
  url       = {https://github.com/yannickwerner/ScenarioReductionForSUC},
}

@Misc{gurobi,
  author  = {{Gurobi Optimization, LLC}},
  title   = {{Gurobi Optimizer}},
  year    = {2024},
  url     = {https://www.gurobi.com},
  version = {11.0.2},
}

@Article{bezanson2017julia,
  author    = {Bezanson, Jeff and Edelman, Alan and Karpinski, Stefan and Shah, Viral B},
  journal   = {SIAM review},
  title     = {Julia: A fresh approach to numerical computing},
  year      = {2017},
  number    = {1},
  pages     = {65--98},
  volume    = {59},
  publisher = {SIAM},
  url       = {https://doi.org/10.1137/141000671},
}

@Article{Lubin2023,
  author    = {Miles Lubin and Oscar Dowson and Joaquim {Dias Garcia} and Joey Huchette and Beno{\^i}t Legat and Juan Pablo Vielma},
  journal   = {Mathematical Programming Computation},
  title     = {{JuMP} 1.0: {R}ecent improvements to a modeling language for mathematical optimization},
  year      = {2023},
  issn      = {1867-2957},
  month     = jun,
  number    = {3},
  pages     = {581–589},
  volume    = {15},
  doi       = {10.1007/s12532-023-00239-3},
  publisher = {Springer Science and Business Media LLC},
}

\end{document}